\documentclass{article}

\usepackage[cmex10]{amsmath}
\usepackage{blindtext}
\usepackage{graphicx}
\usepackage{cite}
\usepackage{rotating}
\usepackage{color}
\usepackage[table,xcdraw]{xcolor}
\usepackage{lscape}
\usepackage{url}
\usepackage[utf8]{inputenc}
\usepackage{cleveref}
\usepackage{eurosym}

\crefformat{equation}{(#2#1#3)}
\crefrangeformat{equation}{(#3#1#4)-(#5#2#6)}
\crefmultiformat{equation}{(#2#1#3)}%
{ and~(#2#1#3)}{, (#2#1#3)}{ and~(#2#1#3)}

\usepackage{nomencl}
\usepackage{etoolbox}
\renewcommand\nomgroup[1]{%
	\ifstrequal{#1}{I}{\item[\textbf{Indices and Sets}]}{%
	\ifstrequal{#1}{P}{\vspace{10pt} \item[\textbf{Parameters}]}{%
	\ifstrequal{#1}{V}{\vspace{10pt} \item[\textbf{Variables}]}{}}}%
}
\date{}

\makeatother

\begin{document}
%
\title{A Three-Stage Stochastic Peer-to-Peer Market Clearing Model with Real-Time Reserve Activation}

\author{
Sigurd Bjarghov\\ Christian \O . Naversen\\ Kasper Thorvaldsen\\ Hossein Farahmand\\
Department of Electric Power Engineering \\
Norwegian University of Science and Technology (NTNU) \\
Trondheim, Norway\\
sigurd.bjarghov@ntnu.no
}
\maketitle

\begin{abstract}
With an increasing share of volatile renewables, more flexibility is needed to balance the electricity system. In order to enable local flexibility, we suggest a three-stage stochastic market clearing model using community-based peer-to-peer (P2P) trading with participation in the short term reserve markets. With this approach, internal intraday (IID) and real-time (IRT) markets are cleared where peers can trade to balance their day-ahead market position under uncertainty in load and activation of reserves. Results show increased participation of 19 and 44 \% in reserve markets when P2P IID and IRT is introduced. We illustrate the application of the model and give advice for the future reserve market design highlighting necessary changes in order to unlock end-user flexibility to the reserve market.
\end{abstract}


\makenomenclature
\setlength{\nomlabelwidth}{1.5cm}


\nomenclature[I]{$P$}{Set of peers, index $p$ and $q$}
\nomenclature[I]{$T_h$}{Set of time steps, index $t$}
\nomenclature[I]{$T_m$}{Set of time steps inside $t$, index $\tau$}
\nomenclature[I]{$S_l$}{Set of 2nd stage load realization scen., index $l$}
\nomenclature[I]{$S_a$}{Set of 3rd stage reserve activation scen., index $a$}


\nomenclature[P]{$W_{pt}^{load}$}{Inflexible load [kWh]}
\nomenclature[P]{$PV_{pt}$}{Photovoltaic production [kWh]}
\nomenclature[P]{$P_{t}^{DA}$}{Day-ahead spot price [$\frac{\text{\euro ct}}{kWh}$]}
\nomenclature[P]{$P_t^{ID}, P_t^{RT}$}{External intraday and real-time price [$\frac{\text{\euro ct}}{kWh}$]}
\nomenclature[P]{$P_{t}^{GT}$}{Grid tariff cost [$\frac{\text{\euro ct}}{kWh}$]}
\nomenclature[P]{$P_{t}^{FCR}$}{Frequency containment reserve price [$\frac{\text{\euro ct}}{kW}$]}
\nomenclature[P]{$P^{P2P,IID}$}{Additional cost of IID P2P trading [$\frac{\text{\euro ct}}{kWh}$]}
\nomenclature[P]{$P^{P2P,IRT}$}{Additional cost of IRT P2P trading [$\frac{\text{\euro ct}}{kWh}$]}
\nomenclature[P]{$\eta_p^{ch},\eta_p^{dis}$}{Battery charging/discharging efficiency [\%]}
\nomenclature[P]{$\Delta \tau$}{Sub-hour time step coefficient [$\frac{\tau}{t}$]}
\nomenclature[P]{$O_p^{min}$}{Minimum battery state of charge [kWh]}
\nomenclature[P]{$O_p^{max}$}{Maximum battery state of charge [kWh]}
\nomenclature[P]{$Q_p^{ch},Q_p^{dis}$}{Maximum battery (dis)charging power [kW]}
\nomenclature[P]{$\pi_{t\tau a}^\uparrow,\pi_{t\tau a}^\downarrow$}{\% of up-/downward reserve bid activated [\%]}
\nomenclature[P]{$\rho_l,\rho_a$}{Probability of scenario $l, a$}

\nomenclature[V]{$\chi_{pt}^{buy},\chi_{pt}^{sell}$}{Purchased/sold DA electricity [kWh]}
\nomenclature[V]{$r_{pt}^\uparrow,r_{pt}^\downarrow$}{Up-/downward DA reserve bid [kW]}
\nomenclature[V]{$\psi_{pqtl}^{sell}$}{IID P2P electricity sold from $p$ to $q$ [kWh]}
\nomenclature[V]{$\psi_{pqtl}^{buy}$}{IID P2P electricity bought from $p$ to $q$ [kWh]}
\nomenclature[V]{$\sigma_{ptl}^{ch},\sigma_{ptl}^{dis}$}{2nd stage battery (dis)charging [kWh]}
\nomenclature[V]{$\sigma_{ptl}^{soc}$}{2nd stage battery state of charge [kWh]}
\nomenclature[V]{$\chi_{ptl}^{ID \uparrow}$}{Intraday balancing electricity bought [kWh]}
\nomenclature[V]{$\chi_{ptl}^{ID \downarrow}$}{Intraday balancing electricity bought [kWh]}

\nomenclature[V]{$\sigma_{pt\tau la}^{ch}$}{3rd stage battery charging [kWh]}
\nomenclature[V]{$\sigma_{pt\tau la}^{dis}$}{3rd stage battery discharging [kWh]}
\nomenclature[V]{$\sigma_{pt\tau la}^{soc}$}{3rd stage battery state of charge [kWh]}
\nomenclature[V]{$\chi_{pt\tau la}^{RT \uparrow}$}{Real-time balancing electricity bought [kWh]}
\nomenclature[V]{$\chi_{pt\tau la}^{RT \downarrow}$}{Real-time balancing electricity bought [kWh]}
\nomenclature[V]{$\psi_{pqt\tau la}^{sell}$}{IRT P2P electricity sold from $p$ to $q$ [kWh]}
\nomenclature[V]{$\psi_{pqt\tau la}^{buy}$}{IRT P2P electricity bought from $p$ to $q$ [kWh]}


\nomenclature[V]{$\delta_{ptl}^{2nd}$}{2nd st. binary variable (dis)charge indicator [0,1]}
\nomenclature[V]{$\delta_{pt\tau la}^{3rd}$}{3rd st. binary variable (dis)charge indicator [0,1]}
\printnomenclature


\section{Introduction} \label{Introduction}

In the last decade, we have witnessed a massive price development in distributed energy resources (DER), both in electricity production and storage. This has lead to a fundamental change in our power systems where production and consumption are widespread in the system. The ever increasing share of decentralized production indicates a need for change in our market design. The adaptive role is pushed towards the demand side to a greater extent, after traditionally being on the production side. Another significant development in the last decade is ICT allowing faster, more accurate and more reliable communications. 

Together, ICT and DER are the two main technological enablers for a decentralized peer-to-peer based market design. Especially peer-to-peer trading and local energy/flexibility markets are emerging as a futuristic, but a possible alternative to the current market structure which has limited capabilities in terms of integrating both large and small scale intermittent renewable electricity production, especially behind the meter. A peer-to-peer market design is quite demanding in terms of resolution. Such a redesign requires an increase in the number of trades/transactions, number of agents, time resolution, and volume resolution. Technological enablers are not the only driver towards decentralized electricity markets. Consumer willingness to participate is key for a consumer-centric market revolution to take place. Engaged consumers are installing PV panels to become less reliant on conventional (often CO$_2$ intensive) power, showing a clear social commitment towards supporting the energy transition and the green shift.

One of the main challenges for the power system when changing from large controllable power plants to numerous small scale intermittent production units is the ability to provide reserve power when demand and supply are at a mismatch. As wind and solar appear to win the profitability race of renewables, demand side needs to take a more significant responsibility of delivering flexibility in the system by participating in reserve markets. Because the demand side mostly consists of many small agents, the number of bids will vastly increase. 

To handle a significant increase in agents, it is natural to assume that there will be decentralized markets where prosumers take an active role in energy trading. This includes interacting with other consumers and prosumers to establish bilateral energy contracts. Potential future decentralized electricity market designs have been described in \cite{Parag2016}, where the authors describe different prosumer era market designs. The main candidates vary from full peer-to-peer markets with complete freedom of trades to organized, community-based markets with a clear set of participants. In addition, hybrid solutions with islanded microgrids or connected microgrids are described. Advantages and disadvantages of the potential market designs are thoroughly described in \cite{Sousa2019}. Full peer-to-peer trading represents the most futuristic and democratic approach but has severe challenges such as high investment costs in ICT for full scalability and reliability under a non-central control scheme. In addition, high computational effort to convergence to a market clearing in a non-pool based market is required, with scalability mentioned as a central issue in \cite{Moret2018} and \cite{Sorin2019}, especially when product differentiation is introduced. As full peer-to-peer market clearing problems often become intractable, community-based peer-to-peer markets with a limited number of peers appear as a more robust approach as scalability is less of an issue. By forming communities or energy collectives \cite{Moret2018a}, a community manager supervises trades and assures convergence to system optimality. Communities can be formed in parts of the grid where there is no congestion. Alternatively, a grid tariff can be enforced if energy collectives choose to trade regardless of grid congestions \cite{Baroche2019}. Common for \cite{Moret2018, Sorin2019, Moret2018a, Baroche2019} is the use of distributed optimization in order to clear the decentralized market, often in the form of Alternating Direction Method of Multipliers (ADMM) or Consensus ADMM (progressive hedging) which converges to a local market energy price given all the agents' utility functions. Another distributed optimization technique known as Consensus + Innovation also shows potential for multiagent coordination \cite{Hug2015}.

There is a number of peer-to-peer pilot projects ongoing such as the Brooklyn Microgrid \cite{Mengelkamp2018}, where excess production from photovoltaics is sold to neighbours in the microgrid in a peer-to-peer marketplace. With decreasing battery costs, a natural next step is to look into the role of battery flexibility peer-to-peer markets. While batteries at peer level reduce overall costs, a community-owned battery reduces import peaks and is more grid friendly \cite{Luth2018}. Similar results are found in \cite{Zepter2019}, although the system cost minimization approach is based on cooperation + savings allocations in the community. Still, intraday trading between peers in order to balance out uncertainty in load and renewable production increases the savings drastically. 

By formulating the peer-to-peer optimization problem as a two-stage stochastic program, \cite{Zepter2019} highlights the need for an internal intraday peer-to-peer market to deal with uncertainty in net load in half-hourly intervals. However, 15-minute market resolutions are quite common in Central Europe and the Nordics are planning to follow in order to deal with high uncertainty and fluctuations in load and renewable production \cite{entsoe_15min}. Real-time markets are becoming more and more relevant to deal with very short term imbalance, and three-stage market clearing models are outperforming two-stage representations of intraday and real-time markets \cite{Abbaspourtorbati2017}, with the three stages representing day-ahead, intraday and real-time markets. Meanwhile, \cite{Hollinger2017} concludes that batteries can increase participation in frequency containment reserve (FCR) markets significantly if the minimum duration of delivery is decreased to 15 minutes from today's 30 minutes, due to the relatively small energy content of batteries. Lack of energy content is also relevant in some cascaded hydropower systems, where reserve bids must respect the available energy content in hydro reservoirs \cite{Chazarra2016}. However, activation of FCR is seldom modelled, although the balancing of the agents' market positions is a complex planning problem when uncertainty in activation is considered.

The contributions in this paper are the following:
\begin{itemize}
    \item We develop a mixed-integer three-stage stochastic programming model to determine optimal bidding in day-ahead and reserve markets for end-users in a community-based P2P market (assuming cooperative scheduling)
    \item We model stochastic activation of reserve market participation from batteries at prosumer level with a 5-minute resolution to analyze end-user capability of delivering reserve in real-time
    \item We model an intraday and real-time peer-to-peer market where peers can trade internally to balance their day-ahead market positions under uncertainty in load and activation of reserve
\end{itemize}

The rest of the paper is organized as follows: \Cref{Market} describes the peer-to-peer and future reserve market design. \Cref{Model} provides a mathematical formulation of the three-stage optimization problem, whereas \Cref{casestudy} describes the performed case study. Results are then presented and discussed in \Cref{discussion} followed by conclusions and future work suggestions in \Cref{conclusion}.
\section{Market design} \label{Market}


\subsection{Community-based peer-to-peer market}


In this paper, a community-based P2P approach is chosen due to the participation in ancillary service markets, in this case, FFR and FCR. End-users are often not aware of their electricity consumption and even less about which grid service markets they should participate in. How to bid, evaluating market positions and dealing with activation of reserves should thus be handled by the community manager on behalf of the end-users. One of the main benefits of this market structure is the reduced need for action from end-users. Both in day-ahead and reserve markets, there is always a chance of placing bids that are sub-optimal, which leads to high penalty costs when the market position is not met. Moving this responsibility to the community manager is a clear advantage and is, to some extent, similar to the aggregator in terms of risk management. 

Inside the community, peers can trade freely with each other. It is assumed that the grid between the end-users is owned by the community. As a consequence, interaction with the external market goes through the community manager who manages the day-ahead and reserve bids of the end-users. Grid tariffs and taxes only need to be paid from import over the community main-meter. In the case of load imbalance, peers can trade to balance their market position in the community internal intraday market (IID). Balance can also be bought on the external intraday market (EID) for a higher price. Finally, peers who have participated in the reserve market get activation signals from the TSO, causing a need to deliver energy in real-time. If not dealt with using their own battery, peers can trade in the internal and external real-time market (IRT and ERT), where the latter being much more expensive. 

When activated, peers need to deliver the energy requested by the TSO, which is random from 0-100 \% of the placed bid. The activation restricts end-users from placing bids in the reserve options market, which they are unable to deliver. One of the tricky consequences of integrating end-users in reserve markets is that they traditionally do not participate in day-ahead spot markets, and thus have no baseline as a reference unlike traditional power plants participating in reserve markets. This is solved by including the end-users in the day-ahead market, letting the community manager bid in the day-ahead market on behalf of the end-user given expected scenarios of net load and activations. 

Note that in this approach, we do not model peer utility functions, meaning that it is, in essence, a community cost minimization optimization problem. Cost/savings allocation is assumed but is not included in this paper. However, prices are set up in a manner which incentivizes self-balancing and trading inside the community in case of imbalance. Implicitly, peers will at first balance their consumption. If this is not possible in a specific scenario, balancing through internal community peer-to-peer trading is the second cheapest option, whereas buying balancing power on the external market is the last and most expensive option.

\subsection{FCR and FFR markets}

It is common practice for the system operator to procure reserve production capacity to secure the stability of the power system. In Europe, the reserve capacity is usually split into several products with different activation times traded in different markets with the system operator as the only buyer. The reserves with the quickest response time, often categorized as primary reserves or (FCR), are automatically activated if there is a change in the system frequency. In Europe, the FCR should be fully activated within 30 seconds of the deviation \cite{balancing_practice}. After this point, the secondary reserves will replace the FCR to free them up for further use. 

There are several different reserve products within the primary reserve category, and these products vary between countries. For instance, the Norwegian TSO Statnett procures symmetric FCR-N reserves for deviations within $\pm$0.1 Hz in addition to upward FCR-D reserves which are used when the frequency drops below 49.9 Hz. Statnett is currently developing an FCR-D product for downward regulation \cite{fcr_d_down}. There is an ongoing effort to harmonize reserve products across Europe to enable efficienct reserve exchange. An example of this is the FCR cooperation between several of the continental European TSOs, where the FCR products are standardized and traded in a single market across countries \cite{fcr_cooperation}.

FCR is commonly traded as a symmetric product, which means that participants have to provide the sold capacity for both up and downward regulation. In relation to the European FCR cooperation project, stakeholders had mixed responses to whether asymmetric FCR products should be offered in the new shared market. The main concern was the additional cost and administrative complexity, which lead to the common FCR market only allowing symmetric products \cite{fcr_cooperation}.

As increasing power fluctuations challenges frequency stability, some TSOs are considering new products falls under the category of synthetic inertia or fast-frequency reserves (FFR). FFR is a system service that delivers a fast power change to mitigate the effect of reduced inertial response so that that frequency stability can be maintained. 
The TSO of Great Britain, National Grid, has launched an FFR-service called "Enhanced Frequency Response" (EFR), and fully contracted reserves must be delivered within one second \cite{Roberts_2018}. 

\begin{figure}[h]
    \centering
    \includegraphics[width=1\columnwidth]{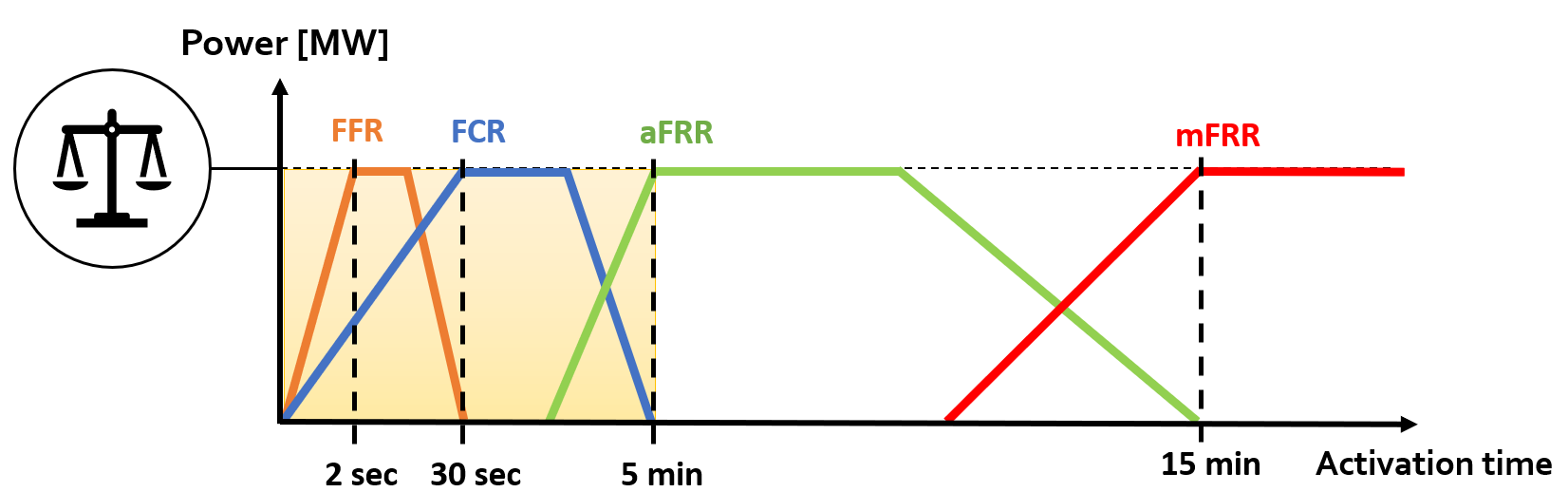}
    \caption{Time horizon for activation of FFR and FCR in Norway. The considered reserve activation time period is shown in yellow.}
    \label{FFR}
\end{figure}

In low-inertia situation in the Nordic system, the traditional hydropower-based FCR-D may not be sufficient to keep the frequency within normal operating range in the event of a dimensioning incident. One alternative solution to maintain secure operation security of power supply is a complimentary service with a response faster than FCR-D. In 2018, the Norwegian TSO, Statnett, conducted a pilot project on FFR where contracted reserves activated within 2 seconds\cite{ffr_pilot} as shown in \Cref{FFR}. The Nordic TSOs are currently considering introducing the FFR product on the ancillary market, and a market release is anticipated in 2020\cite{ffr_pilot}.

In this paper, we assume 5 minute activation periods, but we do not differentiate between FFR and FCR. However, the initial FFR pilot \cite{ffr_pilot} conclude that batteries are well suited to deliver such reserves due to the fast response time.

\subsection{Future reserve market design}
FCR markets have traditionally been symmetric as the service is designed to deal with small, short term deviations in frequency which naturally happen both ways. Symmetric bidding also reduces administration costs for the TSO. However, it is questionable whether or not this is the way to go when designing the future FCR market. Firstly, more and more small scale agents need to participate in order to deal with increased volatility in production due to the increasing renewable share in the EU. There is clearly a vast potential in demand response which needs to be unlocked. However, many of the potential agents have assets which only can participate upward or downward (mostly upward by curtailing load). A significant portion of the potential flexibility will probably not be able to participate if symmetric bids are required. Batteries will also participate less in reserve markets while self-balancing, as a full or empty battery prohibits participation. In addition, minimum bid sizes must be dramatically reduced in order for demand response to participate. This can partially be solved by using communities or aggregators to avoid tiny bids. Still, the minimum capacity of 1 MW which is common in many reserve markets today is quite far away from including any type of residential demand response. Reducing bid size requirements will increase the number of bids and administration costs for the TSO. One could argue that this reduces market efficiency. Moreover, an increased share of renewables indicates that more reserves must be unlocked to meet future reserve requirements.

\section{Model} \label{Model}

The suggested approach is a three-stage stochastic mixed-integer linear programming model, designed to simulate P2P market trading the day-ahead, IID and IRT market. The prosumer needs to take action in three stages; day-ahead, intraday and real-time. A future local market requires more and faster transactions due to deviations in renewables, volatile loads and activation of reserve when participating in external flexibility markets (in this case FFR/FCR). A third stage captures the real-time transactions made to balance these features. Depending on different realizations in load and reserve bid activation uncertainty, first stage decisions on purchased electricity in the day-ahead market and reserve bids needs to be done considering all scenarios.

\begin{figure}[ht]
    \centering
    \includegraphics[width=1\columnwidth]{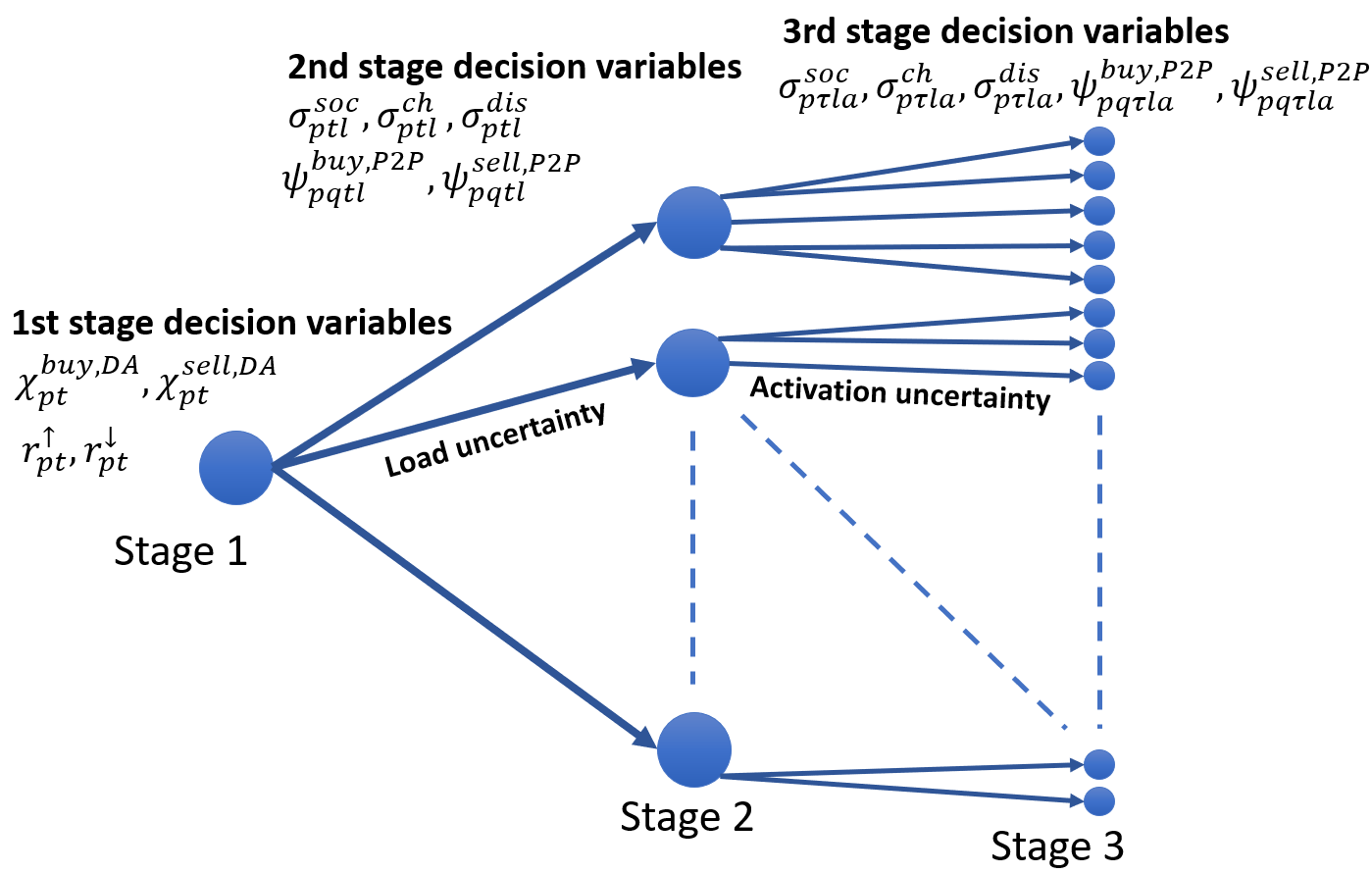}
    \caption{Scenario tree for the three-stage day-ahead, intraday and real-time clearing market.}
    \label{3rdstagebalance}
\end{figure}

The objective function is shown in \Cref{objfun}. We minimize system costs by minimizing the cost of purchased electricity in the day-ahead, EID and ERT markets. Revenue from reserve market participation is subtracted. If activated, a 10\% discount/premium is added to the activated reserve. A cost of P2P trading is added to trades done in the IID and IRT markets.

\begin{align}
    & min \sum_{pt}\chi_{pt}^{buy}(P_t^{DA}+P_t^{GT})-\chi_{pt}^{sell}P_t^{DA} - P_t^{FCR}(r_{pt}^{\uparrow}+r_{pt}^{\downarrow}) \nonumber \\
    &+ \sum_{pta}\Delta\tau\rho_{a}P_t^{DA}(0.9\cdot \pi_t^{\downarrow}r_{pt}^{\downarrow}-1.1\cdot \pi_t^{\uparrow}r_{pt}^{\uparrow}) \nonumber\\
    &+ \sum_{ptl}\rho_l (|\epsilon_{ptl}^{ID}|P_t^{ID} + \psi_{pqtl}^{buy}P_t^{P2P,ID}) \nonumber\\
    &+ \sum_{pt\tau la}\rho_l\rho_a (|\epsilon_{pt\tau la}^{RT}|P_t^{RT} + \psi_{pqt\tau la}^{buy}P_t^{P2P,RT}) \label{objfun}
\end{align}

Set indices are exclusively related to their sets as shown in the nomenclature with $a \in S_a, l \in S_l, \tau \in T_m, t \in T_h$, except for the prosumer indices $p$, where only peers with batteries/PV have battery/PV related variables and parameters. In addition, $q$ is used as the second peer index for P2P trading from peer $p$ to $q$. 
For all realizations of the load uncertainty scenarios $l \in S_l$, the energy balance must be satisfied, linking the two stages as market position represented by the first stage variables $\chi_{pt}^{buy}, \chi_{pt}^{sell}$, $r_{pt}^\uparrow$ and $r_{pt}^\downarrow$ are the same for all scenarios. Depending on the load realization, the prosumer needs to balance his net load to match the day-ahead position. This can be done either by:
\begin{itemize}
    \item Charging or discharging the battery,
    \item Buying or selling electricity on the local P2P intraday market,
    \item Paying an imbalance cost to the external intraday market,
\end{itemize}

which is shown in \cref{eq:ebalancess}. P2P electribity is defined as bought to peer $p$ from peer $q$ \cref{eq:p2p}. 

\begin{align}
    & \chi_{pt}^{buy} - \chi_{pt}^{sell} = W_{ptl}^{load} - PV_{ptl} + \sigma_{ptl}^{ch} - \sigma_{ptl}^{dis} \label{eq:ebalancess} \\
    & + \sum_{q \in P_b} (\psi_{pqtl}^{sell} - \psi_{pqtl}^{buy}) + \chi_{ptl}^{ID\uparrow} - \chi_{ptl}^{ID\downarrow} \hspace{40pt} \forall l, t, p \nonumber\\
    & \psi_{pqtl}^{buy} = \psi_{qptl}^{sell} \hspace{105pt}\forall l, t, q, p, p\neq q \label{eq:p2p}
\end{align}

Reserve market related constraints \crefrange{eq:symmetricreserve}{eq:maxdown} are only defined for prosumers $p \in P_b$, as flexibility (in this case a battery) is needed to provide reserve. Reserve bids are required to be symmetric as described in \cref{eq:symmetricreserve}. The actual activation of bids in the option market requires constraints limiting the reserve bid volume not to exceed available flexibility in time step $t$. The available flexibility has to be measured in two dimensions; energy \crefrange{eq:maxupreserve}{eq:maxdownreserve} and power \crefrange{eq:maxup}{eq:maxdown}. \crefrange{eq:maxupreserve}{eq:maxdownreserve} restrict reserve bids in $t$ to not exceed the available up- and downward state of charge planned for $t+1$, given 100 \% activation in the respective direction. \crefrange{eq:maxup}{eq:maxdown} similarly restrict reserve bids to not exceed available (dis)charging power in $t$. Note that if a charging action is planned, not charging is equal to delivering upward reserve, allowing reserve bids to exceed $Q_p^{ch}$ due to change in the first stage scheduling. 

\begin{subequations}
	\begin{align}
    & r_{pt}^\uparrow = r_{pt}^\downarrow \hspace{168pt} \forall t, p \label{eq:symmetricreserve} \\
    & r_{pt}^\uparrow \leq \sigma_{p(t+1)s}^{soc}\hspace{139pt} \forall l, t, p\label{eq:maxupreserve} \\
    & r_{pt}^\downarrow \leq O_p^{max} - \sigma_{p(t+1)s}^{soc}\hspace{103pt} \forall l, t, p\label{eq:maxdownreserve} \\
    & r_{pt}^{\uparrow} \leq Q_p^{dis} - \sigma_{ptl}^{dis} + \sigma_{ptl}^{ch} \hspace{96pt} \forall l, t, p\label{eq:maxup} \\
    & r_{pt}^{\downarrow} \leq Q_p^{ch} + \sigma_{ptl}^{dis} - \sigma_{ptl}^{ch} \hspace{98pt} \forall l, t, p\label{eq:maxdown}
	\end{align}
\end{subequations}

It is quite common to ignore modelling of activation of FCR/FFR mainly for two reasons; the energy delivered over the time period often has a net value of zero, and conventional power plants can often deliver the requested energy as long as the capacity has been reserved. This approach does not work well for batteries, because the energy actually has to be available, and one can not assume "unlimited" available fuel. When the TSO activates the capacity bid in the same direction for multiple time periods, the energy needs to be available in the battery, also in extreme cases\footnote{Of course a penalty cost can be taken if the activation is not respected, but could lead to exclusion from the market if a trend is spotted by the TSO.}.

The battery is modelled as an electricity storage system with min. and max. state of charge (SOC), (dis)charging power and constant efficiency from the inverter and electrical components \crefrange{eq:soc}{eq:dismax}. The battery starts with a SOC of 0 \%. 

\begin{subequations}
	\begin{align}
    &\sigma_{p(t+1)s}^{soc} = \sigma_{ptl}^{soc} + \sigma_{ptl}^{ch} \eta_{p}^{ch} - \frac{\sigma_{ptl}^{dis}}{\eta_{p}^{dis}}\hspace{60pt}\forall l, t ,p \label{eq:soc}\\
    &O_{p}^{min} \leq \sigma_{ptl}^{SOC} \leq O_{p}^{max}\hspace{100pt} \forall l, t, p \label{eq:socminmax}\\
    &\sigma_{ptl}^{ch} \leq Q_{p}^{ch} \delta_{ptl}^{2nd} \hspace{135pt} \forall l, t, p\label{eq:chmax}\\
    &\sigma_{ptl}^{dis} \leq Q_{p}^{dis}(1-\delta_{ptl}^{2nd}) \hspace{105pt}\forall l, t, p\label{eq:dismax}
	\end{align}
\end{subequations}

In the third stage, reserve market option bids are activated from 0-100\% in scenarios $a \in S_a$. Hence, in all scenarios $a$, the activated reserve must be delivered using the three options:
\begin{itemize}
    \item Adjusting the second stage battery schedule,
    \item Buying or selling electricity in the IRT market,
    \item Buying electricity in the ERT market.
\end{itemize}

The third stage energy balance is shown in \cref{eq:ebalancets} with P2P trading between $p$ and $q$ (battery peers only) in shown in \cref{eq:p2pts}.

\begin{align}
    &\sigma_{pt\tau a}^{dis} - \sigma_{pt\tau a}^{ch} + \pi_{t\tau a}^\uparrow r_{pt}^\uparrow - \pi_{t\tau a}^\downarrow r_{pt}^\downarrow = \sigma_{pt\tau a}^{dis}  - \sigma_{pt\tau a}^{ch} + \chi_{pt\tau a}^{RT\uparrow} \nonumber \\
    &- \chi_{pt\tau a}^{RT\downarrow}+ \sum_{q \in P}(\psi_{pqt\tau a}^{sell} - \psi_{pqt\tau a}^{buy})\hspace{56pt}\forall a, \tau, t, p \label{eq:ebalancets} \\
    & \psi_{pqt\tau la}^{buy} = \psi_{qpt\tau la}^{sell}\hspace{75pt}\forall a, l, \tau, t, p, q, p\neq q \label{eq:p2pts}
\end{align}

The battery is modeled equally in the third stage as in the second stage, with the addition of a new time index $\tau$ representing the five minute intervals of each activation scenario $a$ \crefrange{eq:socts}{eq:dismaxts}. 

\begin{subequations}
    \begin{align}
        &\sigma_{pt(\tau+1)a}^{soc} = \sigma_{pt\tau a}^{soc} + \Delta \tau (\sigma_{pt\tau a}^{ch} \eta_p^{ch} - \frac{\sigma_{pt\tau a}^{dis}}{\eta_p^{dis}}) \hspace{5pt} \forall a, \tau, t, p  \label{eq:socts} \\
        & O_p^{min} \leq \sigma_{pt\tau a}^{soc} \leq O_p^{max} \hspace{90pt}\forall a, \tau, t,  p\label{eq:socminmaxts} \\
        & \sigma_{pt\tau a}^{ch} \leq Q_p^{ch}\delta_{pt\tau la}^{3rd} \hspace{110pt} \forall a,\tau, t, p\label{eq:chmaxts} \\
        & \sigma_{pt\tau a}^{dis} \leq Q_p^{dis}(1-\delta_{pt\tau la}^{3rd}) \hspace{82pt}\forall a, \tau, t, p \label{eq:dismaxts}
    \end{align}
\end{subequations}
\section{Case study} \label{casestudy}

\begin{figure}[ht]
    \centering
    \includegraphics[width=1\columnwidth]{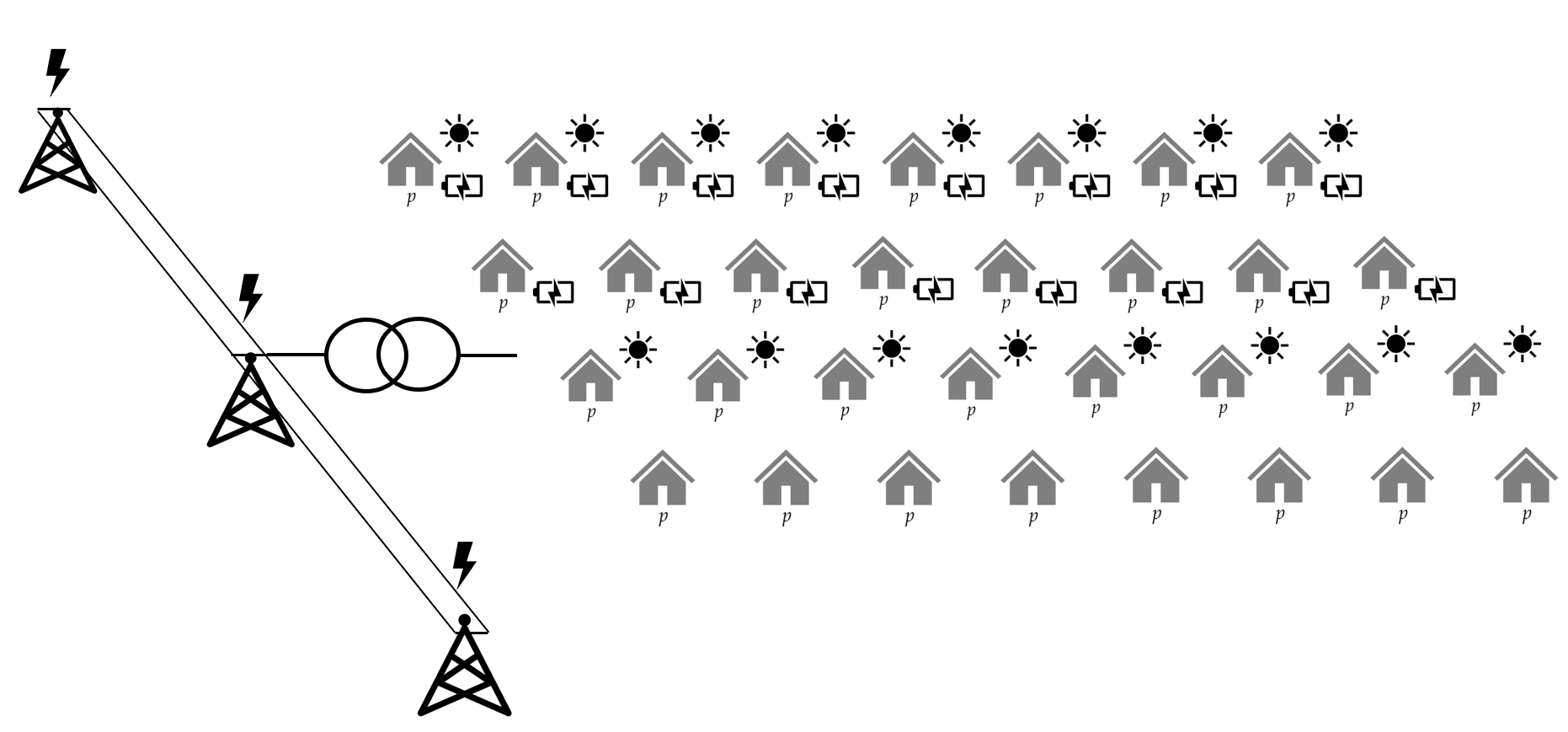}
    \caption{Community setup. All peers can trade in the second stage, whereas only peers with batteries can trade in the third stage.}
    \label{casestudyfigure}
\end{figure}

In order to illustrate, a case study with 32 peers is tested. The peers are a mix of prosumers with batteries and photovoltaic production and consumers with assets distribution shown in \Cref{casestudysetup} and \Cref{casestudyfigure}. The case study is set up in order to analyse the synergy between PV production, battery flexibility and pure consumption. All batteries have 10 kWh with 10 kW maximum (dis)charging power, whereas all rooftop PV is set to 3 kWp.

\begin{table}[h]
\centering
\caption{Share of assets between the peers.}
\begin{tabular}{lcccc}
\hline
\textbf{Share}  & 25\%    & 25\%     & 25\% & 25\% \\
\textbf{Assets} & Battery & Battery + PV & PV   & -    \\ \hline
\end{tabular}
\label{casestudysetup}
\end{table}

We perform three case studies:
\begin{enumerate}
    \item Market clearing with P2P trading in the IID and IRT markets. All peers can trade in the second stage in order to balance under load uncertainty. Only peers with batteries participate in reserve markets, and thus only peers with batteries can trade in the third stage for balancing after realization of activation.
    \item Market clearing with P2P trading only in the IID market. All peers can trade in the second stage to balance under load uncertainty. No P2P trading in the IRT market.
    \item Market clearing without P2P trading. End-users can only use their batteries to balance under load and activation uncertainty. 
\end{enumerate}

In the first stage, all peers make day-ahead bids in the spot market to cover their load. In the second stage, the peers can balance their market position depending on the realization of their net load by P2P trade or battery if available. However, only peers with batteries do P2P trading and balancing in the third stage, as activation only occurs for peers with batteries.

The simulation is done for one day with 24 timesteps of one hour in the first and second stage, whereas there are 12 steps in the third stage (each representing 5 minutes of activation), making a total of 288 time periods. Load data is taken from real households located in Steinkjer, Norway. To generate load scenarios, we use load data from five different days with similar temperatures which represent the forecasted load. All days are taken from the same day of the week but in different weeks. PV is considered deterministic and is generated using irradiation data close to the residences.

Activation scenarios are generated using a normal distribution with an expected value $\mu = 0$ with a standard deviation of $\sigma = 0.3$. All randomly drawn scenarios higher than 100 \% activation are clipped to 100 \%.

Prices for day-ahead spot market and FCR are acquired from NordPool in region NO3 with the corresponding grid tariffs from the relevant DSO. In the case of symmetric bids, FCR prices are taken directly from NordPool. When activated, peers are getting a 10 \% discount and premium for upward and downward reserve energy delivered, respectively. P2P trading is penalized by 5 and 10 \%  of the day-ahead price in the IID and IRT markets, respectively. This is to capture costs related to transactions and to incentivize balancing in the first and second stage, rather than the third stage.
\section{Results and discussion} \label{discussion}

\begin{figure}[ht]
    \centering
    \includegraphics[width=1\columnwidth]{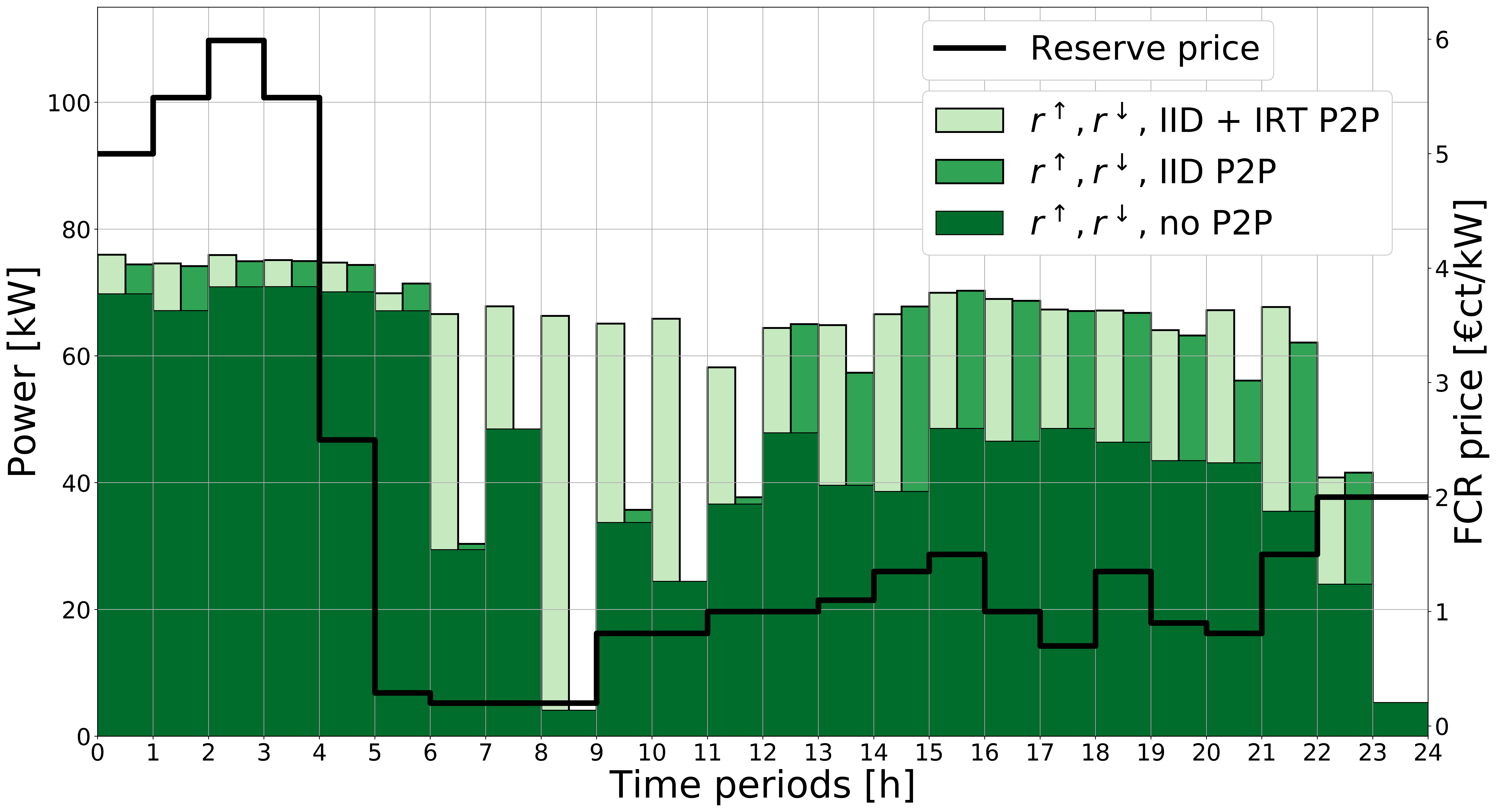}
    \caption{Reserve bids with no, IID only and IID + IRT peer-to-peer trading.}
    \label{reservebids}
\end{figure}

Results show that participation in reserve markets is significantly higher with peer-to-peer trading enabled. From \Cref{reservebids} we see how market participation is similar during the high price periods in the night, utilizing close to all available battery capacity. When FCR prices decrease at 5 am, bids in the no and IID peer-to-peer simulation drop because peers are not able to balance only utilizing their batteries and need to take too high penalty costs if they bid. Reserve bid volumes increase somewhat around noon but stay stable below the peer-to-peer case study where peers are able to bid close to the max of available battery capacity in all hours until 10 pm. The drop in reserve bids is a result of it being more profitable to discharge the battery than to participate in the reserve market, and does not happen if the simulation time horizon does not end. 

When the last two hours are excluded from the results, internal peer-to-peer markets increase participation in reserve markets by 19 and 44 \% in the cases of adding IID and IID + IRT markets, respectively. \Cref{reservebidstable} show an increase in 384 and 883 kW of total bids in the simulated time horizon. The results show that peer-to-peer trading is a powerful market design that significantly increases end-users possibility to participate in reserve markets, due to the ability to meet the market position through balancing. Because some peers have low-load scenarios and others have high, it makes more sense to trade internally than to face penalties from the external market. This kind of synergy is key to unlock the potential of reserve provision using local electricity markets to the system's advantage.

\begin{table}[ht]
\caption{Reserve bid volumes for the three case studies.}
\begin{tabular}{lccc}
\hline
                                                                          & \textbf{No P2P} & \textbf{IID only} & \textbf{IID + IRT} \\
\textbf{\begin{tabular}[c]{@{}l@{}}Reserve \\volume {[}kW{]}\end{tabular}} & 1991            & 2375 (+19\%)             & 2874 (+44\%)              \\ \hline
\end{tabular}
\label{reservebidstable}
\end{table}

\begin{figure}[ht]
    \centering
    \includegraphics[width=1\columnwidth]{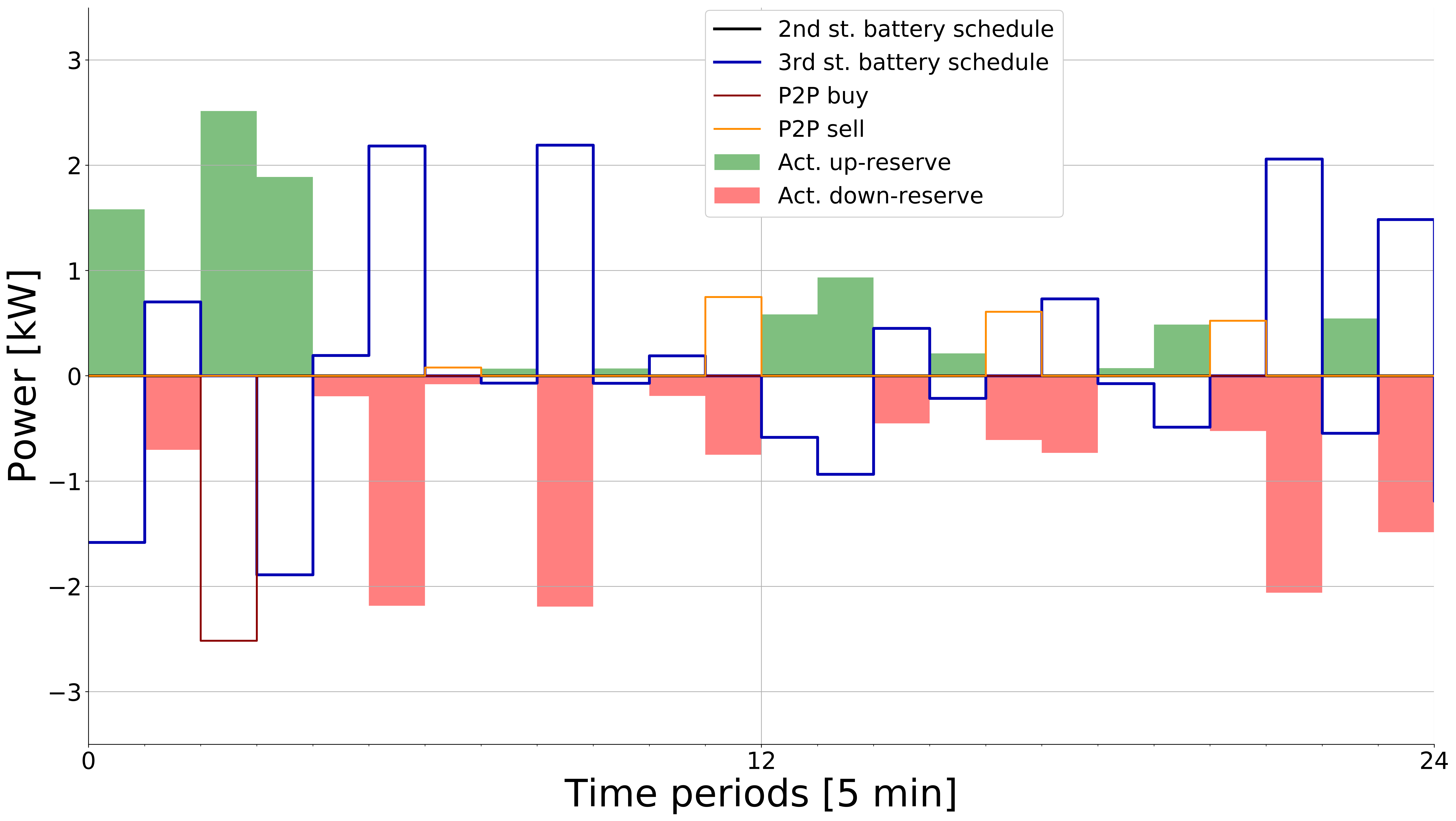}
    \caption{Prosumer reserve activation balancing in the 3rd stage with peer-to-peer trading.}
    \label{TSbalancep2p}
\end{figure}

We visualize the results by showing how the third stage real time market balance is done when reserves are activated. Here, reserve bids are activated by the TSO and peers must deliver the energy required. In \Cref{TSbalancep2p}, an arbitrary example is shown where reserves are activated with different amplitudes up- and downward. Balancing is mostly done using the battery. However, in the third time period, all of the up-reserve is delivered by buying electricity from other peers on the IRT market. In period 12, 17 and 21, the downward activation is met by increasing consumption (selling electricity to other peers). 

\begin{figure}[ht]
    \centering
    \includegraphics[width=1\columnwidth]{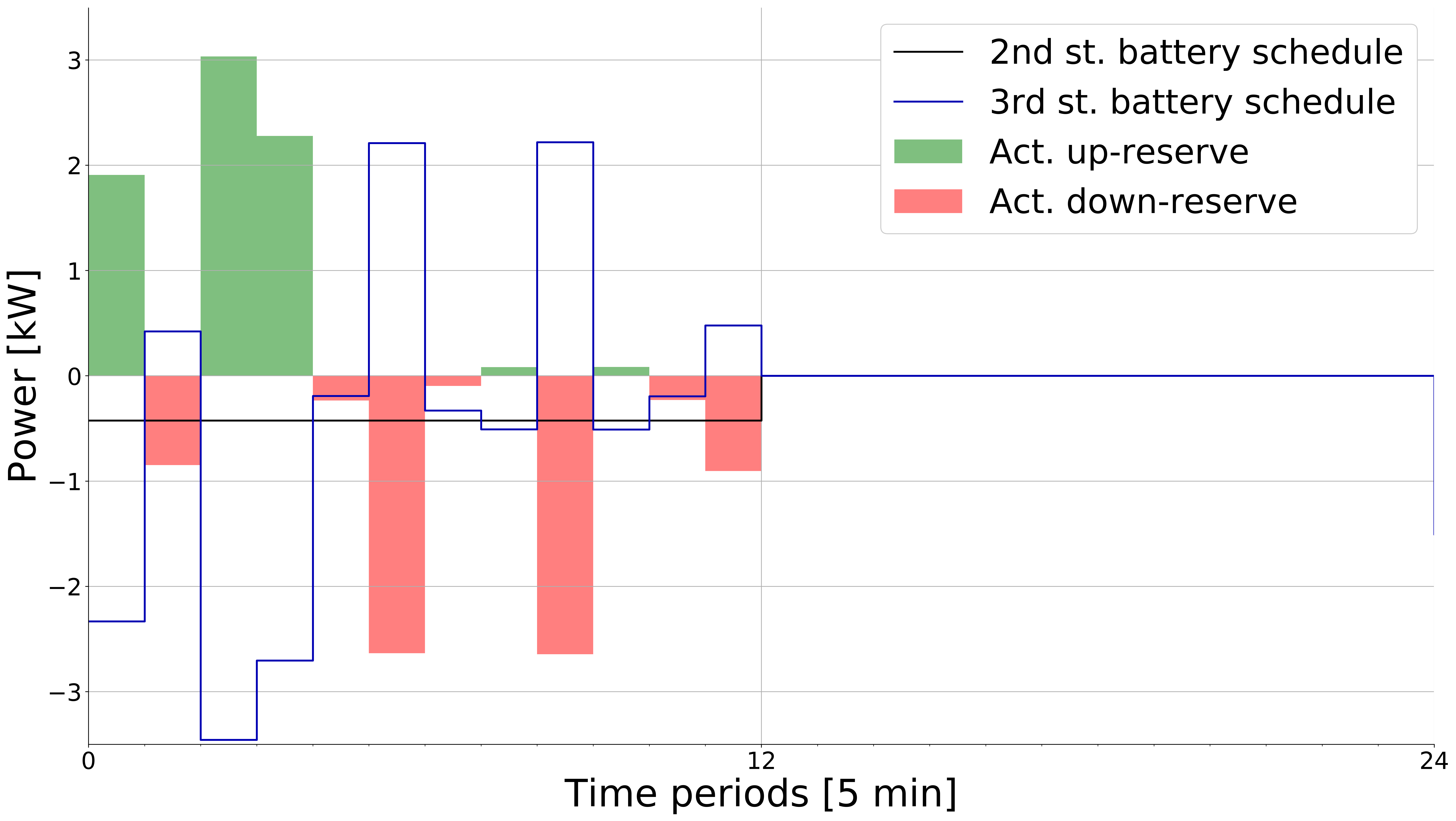}
    \caption{Prosumer reserve activation balancing in the 3rd stage without peer-to-peer trading. No bids are made in the second hour.}
    \label{TSbalancenop2p}
\end{figure}

In the case without peer-to-peer trading, activation can only be met by using the battery or by purchasing the electricity in the external real-time market (for a high cost). In \Cref{TSbalancenop2p}, balance is met only using the battery. Because the battery is already being discharged as part of the intraday scheduling (represented by the black line), the battery has to deviate from that schedule and deliver what is activated. For example, the first period has an activation of $\sim$1.9 kW. Because the battery is already discharging $\sim$0.4 kW to meet the intraday market position, $\sim$1.9 kW has to be added to this by increasing the discharging to $\sim$2.3 kW. In the visualized hour, the net activations are quite close to zero. In the second hour, the net activation is skewed towards downward delivery. Because the FCR price in this hour is quite low, no bid is made due to too high penalty costs in the scenarios where reserves cannot be delivered, unlike in the peer-to-peer case in \Cref{TSbalancep2p}, where some of the reserves are met using peer-to-peer trading.

In order to incentivize end-users to participate in reserve markets, we find that the minimum required bid size must be reduced. A community approach could work as a market player to aggregate flexibility in order to assure acceptable bid sizes. As batteries and demand response outperform traditional power plants, end-user flexibility strike as a strong candidate for flexibility provision in the FFR market.
\section{Conclusion} \label{conclusion}

In this paper, we suggest a three stage community-based peer-to-peer market clearing model in order to capture the volatility that follows when participating in FCR markets, simulating a day-ahead, intraday and real-time market in the first, second and third stage, respectively. We support the idea of a community manager which handles the day-ahead spot market and reserve bids, based on historical data for scenario generation and available flexibility. The community manager also coordinates peer-to-peer transactions when reserves are activated by the TSO. A clear strength in this approach is to move risk from the end-user to a professional agent. 

By considering uncertainty in net load and reserve activation, community-based peer-to-peer trading is used to balance the first stage market position on the day-ahead and reserve market. Results show that reserve bids increase by 19 \% and 44 \% when IID and IID+IRT markets are introduced, respectively. Real-time peer-to-peer markets can work as an efficient tool to enable more end-user participation in FCR.

Future work is to expand the model to include market clearing with competition between peers, using equilibrium models or other decomposition techniques mentioned in the introduction in order to get a clear view of how behaviour impacts the results compared to our cooperation and cost allocation based approach. Furthermore, the case study can be expanded by adding multiple flexibility assets such as flexible demand in heat and electric vehicle charging. In order to incorporate these types of flexibility, asymmetric bids could enable "one-way" flexibility to participate in FCR.

\bibliographystyle{IEEEtran}
\bibliography{IEEEabrv,citations}

\end{document}